\documentclass{article}
\usepackage[utf8]{inputenc}
\usepackage[english]{babel}
\usepackage{csquotes}
\usepackage{biblatex}
\usepackage{hyperref}

\usepackage{biblatex}

\usepackage{ dsfont }
\usepackage{amssymb}
\usepackage{amsmath}        
\usepackage{amsfonts}       
\usepackage{amsthm}         
\usepackage{bbding}         
\usepackage{bm}             
\usepackage{graphicx}       
\usepackage{fancyvrb}       
\usepackage{indentfirst}    
\usepackage{icomma}         
\usepackage{dcolumn}        
\usepackage{booktabs}       
\usepackage{paralist}       
\usepackage{xcolor}         
\usepackage{ dsfont }

\input xy                         
\xyoption{all}
\usepackage{amsmath,amscd}
\usepackage{textcomp}
\theoremstyle{plain}
\newtheorem{veta}{Věta}
\newtheorem{Thm}[veta]{Theorem}
\newtheorem{Prop}[veta]{Proposition}
\newtheorem{Ex}[veta]{Example}
\newtheorem{Cor}[veta]{Corollary}
\newtheorem{Lemma}[veta]{Lemma}

\newtheorem{Remark}[veta]{Remark}

\theoremstyle{plain}
\newtheorem{Def}[veta]{Definition}

\theoremstyle{remark}

\theoremstyle{plain}

\newcommand{\R}{\mathbb{R}}

\newcommand{\N}{\mathbb{N}}

\newenvironment{dukaz}{
  \par\smallskip\noindent
  \textit{Proof}.
}{

\rightline{$\qedsymbol$}
}

\usepackage{blindtext}
\title{Möbius function and modules with finitely many submodules}
\author{Dominik Krasula}

\begin{document}
 \maketitle 
\begin{abstract} 
The article studies the Möbius function for modules and related questions about the finiteness of a lattice of submodules. 

For a module with finitely many submodules, its Möbius function is a function from the lattice of submodules to integers such that for any nonzero submodule, the sum of values over all of its submodules is zero with the value of the zero module being set to one. It is derived from the Möbius function on posets used in enumerative combinatorics, allowing Möbius inversion for functions on lattices of submodules. For a nontrivial and non-semisimple module, its value is zero.  We show how to calculate the Möbius function for a semisimple module based on the cardinalities of endomorphism rings of simple submodules. 

Modules with finitely many submodules are characterised in terms of their semisimple subfactors. It then follows that a distributive module of finite length has finitely many submodules. Rings whose modules of finite length all have finitely many submodules are characterised as rings whose simple modules all have a finite endomorphism ring. Rings with no such simple modules are characterised as rings over which any module with finitely many submodules is distributive.

General results are applied on the bounded path algebras of acyclic quivers over an infinite field. Representations with finitely many subrepresentations then correspond to thin representations. 
\end{abstract}

\smallskip


\smallskip

\noindent \textbf{Author:} RNDr. Dominik Krasula

Charles University, Faculty of Mathematics and Physics, Department of Algebra,

Prague, Czech Republic

krasula@karlin.mff.cuni.cz

ORCID: 0000-0002-1021-7364

\smallskip 

\noindent This work is a part of project SVV-2023-260721

\noindent This research was supported by the grant GA ČR 23-05148S from the Czech
Science Foundation

\section{Introduction}

The Möbius function and the Möbius inversion formula are classical tools in number theory. Later, the Möbius function was defined for arbitrary locally finite posets.  It is the two-sided inverse of the zeta function in the incidence algebra of a given poset. See Section \ref{Definition} for a definition and [12, Prop. 1] for a recursive formula. This combinatorial view of Möbius's function is traditionally associated with G. C. Rota's article [12], ca. 1964. The  \textit{number-theoretic} Möbius function can be interpreted as a particular case for the poset of natural numbers ordered by the divisibility relation. The value of a natural number $n$ corresponds to the value of any couple of numbers whose quotient is $n$. See [12, Ex. 1] for details.

The Möbius inversion formula was first formulated in the mid-1930s in independent works by P. Hall and L. Weisner,  motivated by the study of $p$-groups. Initially, the Möbius function was defined in terms of spanning subsets of a cross-cut in the lattice of subgroups. This definition is equivalent to Rota's definition for a locally finite lattice (called \textit{hierarchy} in Weisner's work); see [12, Thm. 3]. We refer the reader to [12] for a bibliography of early uses of the Möbius function. 

 Since the late 1970s, several generalisations of the Möbius function for categories have appeared, each assuming some \textit{finitness}  condition on the category; see [11] for a bibliography on the subject and detailed discussion of several versions of the Möbius function for categories.  

The category of modules, even when restricted to a full subcategory of modules with finitely many submodules, is too \textit{large} for categorical techniques of Möbius inversion to work. Instead, for a module, we consider the lattice of its submodules, using Rota's definition of the Möbius function; see Section \ref{Definition}. 

This approach coincides with the Möbius function for groups and can be used for any variety of algebras. The case of modules has particularly convenient properties, as any submodule is a kernel of some homomorphism. Thus, as with the classical Möbius function, rather than seeing the Möbius function acting on an ordered pair, a module and its submodule, it can be interpreted as acting on a single module — the corresponding factor module.  Furthermore, the fact that modules with no nontrivial factors have no nontrivial submodules significantly simplifies the calculation of the Möbius function.

However, the introduction of the Möbius function for modules is not motivated by its theoretical properties.  In [7], T.  Honold and A. A. Nechaev used Möbius inversion in algebraic coding theory. They also gave an explicit formula for the Möbius function of finite modules over finite rings.  Honold later used their results in [6] to give a combinatorial characterisation of finite Frobenius rings. M. Greferath and S. E. Schmidt in [5] used Möbius inversion to give a new proof of the MacWilliams theorem for finite Frobenius rings, initially proved in [14] by J. A. Wood. The problem of characterising (one-sided) artinian rings satisfying MacWilliams theorems was solved in  [9]. However, the methods used still need to discuss the finite case separately; see [9, Thm. 2.1].

\medskip

When attempting to generalise the explicit formula, [7, Prop. 2], to representations of fin. dim. algebras, we realised that a similar formula could be obtained over arbitrary rings, assuming that the value is well-defined, i.e., the given module has finitely many submodules; see Section \ref{SecCalculation}. We use the same combinatorial properties as in [7], formulated in Subsection \ref{SecComb}.  However, the proof of formula from [7] depends on the finiteness of a ring. 

\textit{Distributive modules}, i.e.,  modules whose lattice of submodules is distributive, are shown to be closely related to modules with finitely many submodules, see Subsection \ref{SecDist}. Furthermore, when calculating the Möbius function for a fixed module, it is useful to decompose it as a direct sum of \textit{unrelated} modules; see section \ref{SecOrthocyclic}. Such decomposition of a module corresponds to the decomposition of its poset of submodules. This problem naturally appeared early on in the study of distributive modules; see [13, Section 1]. These findings were rediscovered recently in [4]. 

\smallskip

The structure of the article is as follows. Section \ref{Prel} recalls notation and some well-known results used throughout the text.

Section \ref{Sec3} presents a series of module-theoretic observations needed to prove the explicit formula for the Möbius function over a general ring.  In particular, Subsection \ref{SecOrthocyclic} gathers existing results on unrelated modules.  Subsection \ref{SecSemisimple} then studies the finite direct powers of a simple module based on its endomorphism ring.

Section \ref{SecFinMany} then characterises when a module of finite length has a finite lattice of submodules in Proposition \ref{Submodule-finite}. Applying this result to distributive modules in Section \ref{SecDist}, it is shown that a distributive module of a finite length has finitely many submodules. Subsection \ref{SecThin} then shows that a representation of a bound cyclic quiver over an infinite field has finitely many subrepresentations if and only if it is a thin representation, Proposition \ref{thin}. 

Section \ref{Definition} defines the Möbius function for modules and formulates combinatorial properties used in the proof of the explicit formula in Subsection \ref{SecCalculation}.

\section{Preliminaries and notation}  \label{Prel}

This section recalls some properties of modules, posets and representations that will be used throughout the text. By a ring, we always mean an \textit{associative ring with unity}, and $R$ always denotes a ring.   All modules are assumed to be left modules; all ideals are assumed to be left ideals - the right version is analogous.

 For  an  element  $a\in  M$ we write  $Ann(a)=\{r\in R\mid ra=0\}$ viewed as a left ideal. For a module $M$ and $a\in M$ and an $R$-homorphsism $\phi\colon M\to N$ there is an inclusion $Ann(a)\subseteq Ann(\phi(a))$. If the morphism  $\phi$ is an isomorphism, then $Ann(a)= Ann(\phi(a))$.

 Schur's lemma is implicitly used in several proofs. For a simple $R$-module $S$,  a nonzero homomorphism whose domain (codomain) is $S$ is a monomorphism (epimorphism). Consequently,  the ring of endomorphisms, $End_R(S)$,  is a division ring, as all nonzero endomorphisms are isomorphisms. 

\subsection{Poset of submodules}

Let $P$ be a poset and $x, y\in P$. The \textit{interval} $[x,y]$ is a subposet of elements $z$ such as $x\leq z \leq y$. A poset $P$ is \textit{locally finite} if all intervals in $P$ are finite.  It is precisely locally finite posets for which the Möbius function, in the sense of [12], is definable; see Section \ref{Definition}.

For an $R$-module $M$,  its poset of submodules is denoted by $\mathcal{L}(M)$.  It is a complete modular lattice bounded by the zero module $0$ and $M$.  In particular, $\mathcal{L}(M)$ is locally finite if and only if it is finite. In that case, module $M$  has a \textit{finite composition length}, i.e.,   the lattice $\mathcal{L}(M)$ satisfies both ascending and descending chain conditions. If such a module is nonzero, it has nonzero socle, $Soc~M$, and maximal submodules. Minimal  (maximal) submodules correspond to \textit{atoms} (\textit{coatoms}) in $\mathcal{L}(M)$. 

By the Correspondance theorem, for a module $M$ and its submodule $N\leq M$, there is a canonical lattice isomorphism between the interval $[N, M]$ in $\mathcal{L}(M)$ and the lattice $\mathcal{L}(M/N)$. Its restriction to maximal modules containing $N$ is a bijection.

Recall that an  $R$-module is \textit{semisimple} iff  $M$ is generated by simple submodules, iff it is a direct sum of simple modules, iff any submodule of $M$ is a direct summand.

\subsection{Representations of bound quivers}\label{SecPrelQuiv}

This subsection formulates terminology and some properties of representations of bound quivers. We refer to  [ASS] for missing terminology. 
Throughout this text, $K$ always denotes a field.

\smallskip

A \textit{quiver} $Q$ is a quadruple $(Q_0, Q_1, s, t)$ where $Q_0$ is a nonempty  \textit{finite} set of \textit{vertices}, $Q_1$ is a finite set of \textit{arrows} and $s$ and $t$ are two maps $Q_1\to Q_0$ mapping an arrow to its \textit{source} and \textit{target}, respectively.  A vertex is called a \textit{sink} if it is not in the image of the map $s$. A quiver is called \textit{acyclic} if it contains no oriented cycles.

For a field $K$, a quiver $Q$, and an \textit{admissible} ideal $I$ in the path algebra $KQ$, we consider the category of finite-dimensional $K$-linear representations of $Q$ bounded by relations in $I$, denoted by  $rep_K(Q, I)$. A \textit{representation} $M\in rep_K(Q,I)$ consists of a collection of finite-dimensional vector spaces $M_a$ for each $a\in Q_0$ and $K$-linear \textit{structural} maps $M_\alpha$ for each $\alpha\in Q_1$ satisfying the relations given by $I$. A \textit{morphism} $\Phi$ between two representations $N$ and $M$ is given by a collection of $K$-linear maps $\phi_a\colon N_a\to M_a$ for each $a\in Q_0$ commuting with \textit{structural maps}. For a \textit{bounded path algebra} $KQ/I$, the category $mod\text- KQ/I$ is equivalent to $rep_K(Q, I)$.

\smallskip 

There is a bijection between the isomorphism classes of simple modules and $Q_0$. We denote $S(a)$ the simple representation such that $S(a)_a= K$ and $S(a)_b=0$ for $a\neq b\in Q_0$. The endomorphism ring of a simple representation is isomorphic to $K$.

\begin{Lemma}[ASS, Lemma III.2.2]\label{ASSSocle} Let $Q$ be a quiver, $I$ an admissible ideal in $KQ$ and let  $M\in rep_K(Q,I)$. Then

    (a) $M$ is semisimple iff $M_\alpha=0$ for all arrows $\alpha\in Q_1$.

    (b) $Soc~M$ is a subrepresentation of $M$ such that $Soc~M_a=M_a$ if $a$ is a sink and \[Soc~M_a=\bigcap_{\substack{\alpha \in Q_0\\ s(\alpha)=a}} Ker(\alpha),\]
otherwise. Structural maps in $Soc~M$ are restrictions of structural maps in $M$.
\end{Lemma}

\section{Poset of submodules} \label{Sec3}

This section gathers the module-theoretic properties used in the sequel.  Results from Subsection \ref{SecOrthocyclic} will allow us to reduce the calculation of the Möbius function to the case of finite direct powers of simple modules.  These modules are then discussed in Subsection  \ref{SecSemisimple}.

\subsection{Unrelated modules}\label{SecOrthocyclic}

The Möbius function is multiplicative; see [12, Prop. 3.5] for combinatorial details or Lemma \ref{productlemma} for a module-theoretic version. Thus, our starting point is to ask when the poset of submodules of a direct sum of modules can be naturally viewed as a product of their respective posets of submodules, as formalised in Definition \ref{Defunrelated}. This question naturally appears in various contexts, so the terminology is not fixed; see [4, Remarks 5.2]. As far as the article's author is aware, W. Stephenson was the first to give a general characterisation of such families of modules.  Following [13], we will call them \textit{unrelated}. The unrelatedness of a family of distributive modules is an immediate sufficient condition for their direct sum to be distributive. It is also necessary, as shown in [13, Prop. 1.3] and reproven in [4, Thm. A].

\begin{Def}\label{Defunrelated}
    Let $\mathcal{M}=(M_i)_{i\in I}$ be a family of modules for an index set $I$.

    We say that modules in $\mathcal{M}$  are \emph{unrelated} if any $L\subseteq \oplus M_i$ is equal to $\oplus (M_i\cap L)$.
\end{Def}
A natural counter-example is a direct sum of two isomorphic modules. An immediate positive example is as follows.
\begin{Ex}\label{indecomposable}
     Let $R_1$ and $R_2$ be two rings and let $M$ be a finite-length $R_1\times R_2$-module. Module $M$ is a direct sum of modules $M_1:=(1_{R_1},0)M$ and $M_2:=(0,1_{R_2})M$, where $1_{R_1}$ and $1_{R_2}$ denotes units in rings $R_1$ and $R_2$ respectively. Modules $M_1$ and $M_2$ are unrelated.  
\end{Ex}
The following theorem discusses existing characterisations of unrelatedness. 
\begin{Thm}\label{ThmUnrelated} Let $\mathcal{M}=(M_i)_{i\in I}$ be a family of modules for an index set $I$.  The following are equivalent:

(1) Modules in $\mathcal{M}$ are unrelated.

(2) For any $i\neq j\in I$, modules $M_i$ and $M_j$ are unrelated.

(3) For any $i\neq j\in I$, there are no nonzero homomorphisms between subfactors of $M_i$ and $M_j$.

(4) For any $i\neq j\in I$, $M_i$ and $M_j$ have no isomorphic nontrivial subfactors.

(5) For any $i\neq j\in I$, $M_i$ and $M_j$ have no isomorphic simple subfactors.

(6) For any $i\neq j\in I$, and any $(x_i, x_j)\in M_i\oplus M_j$, the left ideals $Ann(x_i)$ and $Ann(x_j)$ are comaximal.
\end{Thm}
\begin{dukaz}
    It is clear that (2) implies (1). The opposite implication follows from the remaining conditions, showing that unrelatedness can be tested pairwise.

Conditions (5), (6), and (1) are equivalent according to [4, Prop. 4.3].

Equivalence of (1) with (4) follows from [13, Prop. 1.2].

It is clear that (3) implies (4) and (4) implies (5).

It remains to show that (4) implies (3). Suppose that (4) holds and there are subfactors $Q_j$ and $Q_j$ of $M_i$ and $M_j$ respectively and $\phi\colon Q_i\to Q_j$ is a homomorphism. Then $Q_i/Ker(\phi)\cong Im(\phi)$ by the homomorphism theorem. Thus by (4), $Im(\phi)$ is the zero module, i.e., $\phi$ is the zero morphism. 
\end{dukaz}

For our purposes, the following special case will be crucial. 
\begin{Cor} \label{Semisimple}
    Let $M$ be a semisimple module,  $S$ be a simple module and $t\geq 1$ a natural number. 

    Then $M$ and $S^t$ are unrelated iff none of the submodules of $M$ is isomorphic to $S$.
\end{Cor}

\subsection{Finite direct powers of simple modules}\label{SecSemisimple}

This section characterises when a semisimple module has only finitely many submodules.  Following Corollary \ref{Semisimple}, it is enough to discuss modules of type $S^t$ for some simple module $S$ and  $t\in \N$. If $S$ is finite, then module $S^t$  has only finitely many modules. In general, the number of submodules of a module of type $S^t$ follows from the cardinality of the division ring $End_R(S)$ as shown in Lemma \ref{S^2} and Corollary \ref{maximal}. The following example shows that an infinite simple module may have a finite ring of endomorphisms. 

\begin{Ex}
    Let $F$ be a finite field and $\kappa$ an infinite cardinal, $n\in \N$.
    
(1) Module $S:=F^{(\kappa)}$ is an infinite  simple  $End_F(S)$-module with a finite ring of endomorphisms.

(2) Module $T:=F^n$ is a simple $End_F(t)$-module and its endomorphisms are in bijection with elements of $F$.
\end{Ex}
\begin{dukaz}  We only prove (1). Let $R:= End_F(S)$ and consider $\phi \in End_R(S)$, i.e., an $R$-linear map $F^{(\kappa)}\to F^{(\kappa)}$ commuting with any endomorphism of $F^{(\kappa)}$. Then, in particular, it commutes with any restrictions of an endomorphism of $F^{(\kappa)}$ to a finite-dimensional subspace. Thus, $\phi$ must be a map that multiplies an element by some scalar $f\in F$.

On the other hand, any such map is an endomorphism, as elements from $R$ commute with multiplication by scalars. Thus $End_R(S)\cong F$.
\end{dukaz}
\begin{Remark}
    If $R$ is a semilocal ring, i.e., its factor by the Jacobson radical is semisimple, then a simple $R\text-$module is finite if and only if its endomorphism ring is. 
\end{Remark}
\begin{Lemma} \label{S^2}
    Let $S$ be a simple $R$-module such that $End_R(S)$ is infinite. 

    Then $\mathcal{L}(S^2)$ is infinite. 
\end{Lemma}
\begin{dukaz}
We fix some nonzero element $a\in S$, and for each nonzero $\phi\in End_R(S)$, we consider a cyclic module $R(a,\phi(a))$. This module is simple: consider a map 
 \[  R\to R(a,\phi(a)) ~~~~~~ r\mapsto (ra,r\phi(a))\]
with  kernel $Ann(a)\cap Ann(\phi(a))$. Because $\phi$ is an isomorphism  we get $Ann(a)=Ann(\phi(a))$, hence $R/Ann(a)\cong R(a,\phi(a))$ is a simple module.

Let  $\phi, \psi\in End_R(S)$ be two isomoprhisms. If $R(a,\phi(a))=R(a,\psi(a))$ then in particular $R(a,\phi(a))$  also contains $(a,\psi(a))$, hence $(0,\phi(a)-\psi(a)) \in R(a,\phi(a))$. Then we get the following  inclusions of modules 
\[0 \subseteq  R(0,\phi(a) -\psi(a)) \subsetneq R(a,\phi(a)). \]
Because $R(a,\phi(a))$ is a simple module we get $R(0,\phi(a) -\psi(a))=0$, so in particular $\phi(a)=\psi(a)$. But $a$ generates simple module $S$, hence $\phi=\psi$.
\end{dukaz}

\smallskip

We can easily calculate the number of submodules of any fixed length for a simple module with a finite endomorphism ring. We start by counting simple submodules. 
\begin{Lemma}\label{CountingSimples}
    Let $S$ be a simple $R$-module $t,q\in N$ such that $|End_R(S)|=q$.

    Then $S^t$ contains $1+q+q^2+\dots+ q^{t-1}$ simple submodules. 
\end{Lemma}
\begin{dukaz} 
Consider the semisimple ring $Q:=End_R(S^t)$ isomorphic to the ring of $t\times t$-matrices over division ring $End_R(S)$. Then, we can view $Hom_R(S, S^t)$ as a \textit{simple} $Q$-module where the action of an element of $End_R(S^t)$ is given by post-composition. Hence $Hom_R(S,S^t)$ is isomorphic (as a $Q$-module) to a $1\times t$-matrix module over the division ring $End_R(S^t)$.
Thus, if $End_R(S)$ is a finite field with $q$ elements, there  are $q^t$ $R$-homomorphisms from $S$ to $S^t$. In particular, $q^t-1$ monomorphism $S\hookrightarrow S^t$. 

Any simple submodule $T$ of $S^t$ is isomorphic with $S$, so there are $q-1$ morphisms in $Hom_R(S,S^t)$ with image $T$ corresponding to non-zero elements in $End_R(S)$. Thus there are 
\[\frac{q^t-1}{q-1}=1+q+q^2+\dots +q^{t-1}\]
distinct, simple submodules of $S^t$.
\end{dukaz}
There is a bijection between the simple and maximal submodules. 
\begin{Cor}\label{maximal}
   Let $S$ be a simple $R$-module $l,t,q\in N$ such that $|End_R(S)|=q$.

   Then $S^t$ contains 
\[\frac{s_{t-l+1}+\dots +s_t}{s_1+\dots+ s_l}\]
   submodules of length $l$, where $s_t=q^{t-1}+q^{t-2}+\dots + q+1$
\end{Cor}
An analogous statement is well known for abelian $p$-groups (here, we use  $q$ instead of $p$); see [2, 48-49]. Once we know the number of simple submodules, the proof is similar. The set of simple submodules $T_1,\dots, T_k$ of $S^t$ is \textit{independent} iff $l(\sum^k_{i=1} T_i)=k$, or equivalently, iff no $T_i$ is in the submodule generated by the remaining modules. Then, we can calculate the number of submodules of fixed length $l$ by calculating the number of \textit{independent} sets of size $l$ contained in $S^t$ and $S^l$.

\begin{Cor} \label{InfSocle}
Let $S$ be a simple module and $t\in \N$ Then $\mathcal{L}(S^t)$ is finite if and only if $End_R(S)$ is finite. 
\end{Cor}

\section{Modules with finitely many submodules} \label{SecFinMany}
For an $R$-module $M$, if $\mathcal{L}(M)$ is finite, then $M$ is a finite length module. If $R$ is finite, the opposite implication is also true. As seen in the previous subsection, over a general ring, e.g., an infinite field, even a finite-length semisimple module can have infinitely many submodules.  

The main result of this section is Proposition \ref{Submodule-finite} characterising when a module has finitely many submodules.  Corollary \ref{CorFinend} then characterises rings whose finite-length modules all have finitely many submodules.   Proposition \ref{Submodule-finite} is used to show that a distributive module of finite length always has finitely many submodules. We also use it to provide a new proof that a representation of a bounded quiver over an infinite field has finitely many subrepresentations iff it is a thin representation. Example \ref{Counterex} shows this is untrue for quivers with oriented cycles, even if the bound quiver algebra is finite-dimensional.

\subsection{Criterion} 

\begin{Prop}\label{Submodule-finite}
    Let $M$ be a finite-length $R$-module such that $\mathcal{L}(M)$ is infinite. 

    Then, there exists a simple $R$-module $S$ and a submodule $K\leq M$ such that  $M/K$ contains a submodule isomorphic to $S^2$, and $End_R(S)$ is infinite.
\end{Prop}
\begin{dukaz}
Consider set
\[\mathcal{M}:=\{M'\leq M \mid \mathcal{L}(M') ~is~infinite \}\]
partially ordered by inclusion. Module   $M$ is artinian, so $\mathcal{M}$ has a minimal module $N$.     By the minimality of $N$, no maximal submodule of $N$ contains infinitely many submodules. Because $N$ has a finite length, this implies that $N$ has infinitely many maximal submodules.

Let $N_0$ be a maximal submodule of $N$. By the minimality of $N$ in $\mathcal{M}$, the poset $\mathcal{L}(N_0)$ is finite, and so is the set
\[\{N_0\cap N'\mid N' ~maximal~in~N\}\subseteq \mathcal{L}(N).\]
 But $N$ has infinitely many maximal submodules, so there exists a submodule $K\leq N'$ such that the set
\[  \mathcal{N}:=\{ N'~maximal~in~ N \mid N'\cap N_0=K\}\]
is infinite. Note that $K$ is maximal in $N_0$ and any module from $\mathcal{N}$.

The canonical projection 
$ \pi\colon M \rightarrow M/K$ then induces a lattice isomorphism between interval $[K,M]$ in $\mathcal{L}(M)$  and  $\mathcal{L}(M/K)$. 
Modules from $\mathcal{N}$  all contain $K$ as a maximal submodule, so their images are distinct, simple modules in $M/K$.
 \end{dukaz}
\begin{Cor} \label{CorFinend}
Let $R$  be a ring. 

All finite-length $R$-modules have finitely many submodules if and only if all simple $R$-modules have a finite ring of $R$-endomorphisms. 
\end{Cor} 
All \textit{finitely generated} $R$-modules have finitely many submodules iff the regular module $R$ is of finite length. This implies that $R$ is left artinian, and the factor $R/rad~R$ is finite. Hence, the ring $R$  is finite. 

We end this section by discussing a particular case of Proposition \ref{Submodule-finite} in the case of noetherian semiperfect rings.
\begin{Remark} \label{RemIov}
    It was recently observed in [8, Cor. 2.3] that over an indecomposable left and right artinian ring, either all simple modules are finite or have the same infinite cardinality. This observation can be generalised for indecomposable semiperfect left and right noetherian rings [10, Thm. 9].     Note that a semiperfect ring can be decomposed as a direct product of finitely many indecomposable rings, sometimes called \textit{blocks} in literature. This decomposition is unique up to isomorphism and reordering.

    Let $R$ be a semiperfect left and right noetherian ring with a block-decomposition $R\cong B_1\times \dots \times B_k$ and $M$ a finite-length $R$-module. Following Example \ref{indecomposable}, we decompose  $M$ as a direct sum of unrelated modules $M_1\oplus\dots\oplus M_k$ such that for each $M_i$, the action of $B_j$ is nontrivial if and only if $i=j$. Furthermore, $R$ is semilocal, so a simple module is finite if and only if its endomorphism ring is.

    Then either all simple $B_i$-modules are finite, and by Corollary \ref{CorFinend}, module $M_i$ has finitely many submodules if and only if it has finite length Or all simple $B_i$-modules are infinite, and by Proposition \ref{Submodule-finite},  $M_i$ has finitely many submodules if and only if it has a finite length and all non-zero epimorphic images of $M_i$ have simple socles.    
    \end{Remark}

\subsection{Distributive modules}\label{SecDist}

Recall that a module is distributive iff its lattice of submodules is distributive. By [3, Thm. 1], a module $M$ is distributive if and only if all of its factors have a square-free socle. In particular, the assumptions of Proposition \ref{Submodule-finite} are satisfied. 
\begin{Cor}\label{CorDis}
 Distributive modules of finite length have only finitely many submodules. 
\end{Cor}
For a large class of rings, the opposite implication is also true. 
\begin{Remark}
    If $R$ is a ring such that all of its simple modules have infinite endomorphism rings, then all modules with finitely many submodules are distributive.

    Following Remark \ref{RemIov}, if $R$ is a semiperfect noetherian ring, a module with finitely many submodules can be written as a direct sum of a finite module and a distributive module. 
\end{Remark}

    \subsection{Representations of bound quivers}\label{SecThin}
   
Regarding Proposition \ref{Submodule-finite}, it seems natural to ask whether it is enough only to investigate socles of \textit{some} factors - such as socles in the socle series.
\begin{Ex}\label{A3}
  Let $K$ be an infinite field and  $Q\colon 1 \xrightarrow{\alpha} 2$ 
    a quiver with  representation $K^2\xrightarrow{(1~0)} K$. 
    
Using Lemma \ref{ASSSocle}, one can see that both $Soc(M)$ and $M/soc(M)$ have finitely many submodules.  On the other hand, $M/S(2)$ has infinitely many submodules, and so does $M$.
    \end{Ex}

We now present the main proposition of this subsection.

\begin{Prop} \label{thin}
    Let $Q$ be an acyclic quiver, $K$ an infinite field,  $I$ an admissible ideal in $KQ$ and let $M\in rep_K(Q,I)$.

    Then $\mathcal{L}(M)$ is finite if and only if $M$ is thin. 
\end{Prop}

\begin{dukaz}
Using Lemma \ref{ASSSocle} and Corollary \ref{InfSocle}, a semisimple representation has finitely many subrepresentations if and only if it is a thin representation. Any factor of a thin representation is thin; thus, its socle has finitely many submodules.

Now assume $M$ is not thin, i.e., there is a vertex $a\in Q_0$ such that $dim_K(M_a)=t>1$. If $a$ is a sink, then by Lemma \ref{ASSSocle}, $Soc~M$ contains an isomorphic copy of $S(a)^t$; thus, it has infinitely many subrepresentations.

Assume $a$ is not a sink and let $b_1,\dots, b_k$ be the set of all vertices that are targets of arrows with the source $a$. To each $b_i$, we assign a quiver $Q^i$, defined as the minimal full subquiver containing $b_i$ such that all arrows with their source in $Q^i$ also have a target in $Q^i$. Because $Q$ is acylic, the minimality of $Q^i$ implies that it does not contain vertex $a$. 

For a quiver $Q'=\cup_{i=1}^k Q^i$ we define a representation  $R_M(Q')$. A  linear space (map) in $R_M(Q')$ is the same as in $M$ if the corresponding vertex (arrow) is in $Q'$ and zero otherwise. Because all arrows of $Q'$ with their source in $Q'$ also have their target in $Q'$, representation $R_M(Q')$ is a subrepresentation of $M$.

We claim that factor representation $M/R_M(Q')$ has a socle with infinitely many subrepresentations. Consider an arrow  $\alpha$ whose source is $a$. The corresponding structural map  in the factor representation $M/R_M(Q')$  is a zero  because its codomain, $(M/R_M(Q'))_{t(\alpha)}$, is the zero $K$-vector space. So, using again Lemma \ref{ASSSocle}, the socle of  $M/R_M(Q')$ contains an isomorphic copy of $S(a)^t$. 
\end{dukaz}
The following example shows that the above proposition does not hold for quivers with oriented cycles even when the given bounded path algebra is finite-dimensional. 
\begin{Ex} \label{Counterex} For infinite filed $K$ and a quiver $Q\colon \xymatrix@=10pt{1 \ar@/^/[r]^\alpha & 2 \ar@/^/[l]^\beta} $ consider representation $M\colon \xymatrix@=10pt{K \ar@/^/[r]^{(1~ 0)} & K^2 \ar@/^/[l]^{(0~1)}}$ of the Frobenius algebra $A:=KQ/I$, where $I$ is an admissible ideal given by a relation $\alpha\beta = 0$.

Representation  $M$  is not thin, but it has only finitely many submodules by Proposition \ref{Submodule-finite}. Its socle is simple, and there is no subrepresentation with dimension vector $[1,0]$, so all non-trivial factors are thin and thus have only finitely many submodules by Proposition \ref{thin}.
\end{Ex}

\section{Möbius function} \label{Definition}

This section recalls the combinatorial definition of the Möbius function and defines its version for modules. It then determines the values for modules with finitely many submodules. The article [12] is used as a reference.

\smallskip 

 For a locally finite poset $P$ and a field $\R$, we consider an \textit{incidence algebra} $\R P$ consisting of all real-valued functions with domain $P^2$. For two elements $\alpha, \beta \in \R P$ the multiplication $*$ is then defined
\[(\alpha*\beta)(x,y)=\sum_{\substack{x\leq z\leq y}}\alpha(x,z)\beta(z,y),\]
if $x\leq y$ and zero otherwise. Kronecker delta is then the multiplicative unit in $\R P$.

The \textit{zeta function} $\zeta_P(x,y)$ is defined as 1 if $x\leq y$ and 0 otherwise. By [12, Prop. 1], the zeta function has a two-sided inverse, called the \textit{Möbius function}, denoted by $\mu_P$.

\begin{Def}
    Let $M$ be an $R$-module such that $\mathcal{L}(M)$ is finite.  
    
    Then $\mu_R(M)$ is an integer defined  recursively by setting  $\mu_R(0):=1$ and if $M\neq 0$ then  $\mu_R(M)$ is the unique integer such that 
     \[\sum_{\substack{N\in \mathcal{L}(M)}} \mu_R(N)=0.\]    
\end{Def}
Note that $\mu_R(-)$ is not a function as the $R$-modules with finitely many submodules do not form a set.  In applications such as Proposition \ref{inversion}, one usually works with one fixed module and its submodules, which then form a set, and our definition then coincides with the definition for posets in the following sense.
\begin{Remark}
 Let $N\leq M$ be two $R$-modules such that $\mathcal{L}(M)$ - and thus also $\mathcal{L}(N)$ - is finite, then 
 \[ \mu_R(N)= \mu_{\mathcal{L}(N)}(0,N) = \mu_{\mathcal{L}(M)}(0,N).\]
The correspondence theorem implies that 
 \[\mu_{\mathcal{L}(M)}(N,M) = \mu_R(M/N).\]
\end{Remark}

\begin{Ex}
    If $S$ is a simple module, then $\mu_R(S)=-1$.  Let $M$ be a module of length 2 with $\mathcal{L}(M)$ finite. If $M$  has a simple socle then $\mu_R(M)=0$. If $M$ is semisimple with $n$ simple submodules then $\mu_R(M)=n-1$.
\end{Ex}

\subsection{Combinatorial properties of Möbius function} \label{SecComb}

The following three lemmas are a module-theoretic reformulation of well-known properties of the Möbius function. The reference for this subsection is [12]. Lemma \ref{non-semisimple} is usually attributed to P. Hall and Lemma \ref{coatom} to L. Weisner.  
 
\begin{Lemma}[12, Prop. 5.2] \label{non-semisimple}
Let $M$ be a nonzero $R$-module. 

If $M$ is not semisimple, then $\mu_R(M)=0$.    
\end{Lemma}
In the following subsection, we will see that the opposite implication is also true. In a general finite lattice $\mathcal{L}$, there might be an element $x$ that is a join of atoms, yet $\mu_\mathcal{L}(0,x)=0$.
\begin{Ex}
   Consider a 6-element lattice $\mathcal{L}$ with the upper and lower bounds 0 and 1, three atoms  $a,b,c$ and two more elements $a\vee b$  and $b\vee c$ with $\mu_\mathcal{L}(0,a\vee b) =\mu_\mathcal{L}(0,b\vee c)=1$. 
   
   In this lattice,  $1=a\vee b\vee c = a\vee c$ and $\mu(0,1)=0$.
\end{Ex}
\begin{Lemma}[12, Prop. 3.5]  \label{productlemma} 
    Let $M, N$ be two unrelated modules such that $\mathcal{L}(M)$ and $\mathcal{L}(N)$ are finite posets.    

    Then, for any $L\leq M\oplus N$ we have 
    \[\mu_R(L)=\mu_R(M\cap L)\cdot \mu_R(N\cap L). \]
\end{Lemma}
\begin{Lemma}[12, Prop. 5.4]\label{coatom}
Let $M$ be an $R$-module such that $\mathcal{L}(M)$ is finite, and let $T$ be a simple submodule.  Then

\[\mu_R(M)= - \sum_{\substack{T\nleq N\leq M\\  N~maximal}} \mu_R(N).\]
\end{Lemma}

\subsection{Calculation of the Möbius function}  \label{SecCalculation}

We assume that all modules in this subsection have only finitely many submodules.  By Lemma \ref{non-semisimple}, the Möbius function of a nonzero module that is not semisimple is zero. 

Let $M$ be a semisimple module. Then there exists a decomposition  \begin{gather}\label{decomposition}
      M= S_1^{t_1}\oplus \dots S_n^{t_n},
  \end{gather}
    such that $S_1,\dots, S_n$ are non-isomorphic simple modules and $t_i$ natural numbers.

    Becuase we asssume that $\mathcal{L}(M)$ is finite, for any $i\leq m$, the inequality  $t_i>1$ implies that $End_R(S_i)$ is finite by Corollary \ref{InfSocle}.

    By Corollary \ref{Semisimple}, we have
\[ \mathcal{L}(M) \cong \mathcal{L}(S_1^{t_1}) \times \dots \times \mathcal{L}(S_t^{t_n}), \]
so by repeatedly applying  Lemma  \ref{productlemma} ,we get 
\[\mu_R(M) =\prod_{\substack{1\leq i\leq n}}\mu_R(S_i^{t_i}).\]

The following lemma then completes the calculation of the Möbius function.
\begin{Lemma}\label{formula}
    Let $S$ be a simple $R$-module, $t$ and $q$ natural numbers such that $|End_R(S)|=q$. Then
\[\mu_R(S^t) = (-1)^{t} q^{\frac{t(t-1)}{2}}.\]
\end{Lemma}
\begin{dukaz}
We apply Lemma \ref{coatom}.  By Corollary \ref{maximal}, we get that $S^t$ contains $1+q+\dots +q^{t-1}$ maximal submodules. 
   
   Let $T$ be some simple submodule of $S^t$. There is a bijection between maximal submodules of $S^t$ containing $T$ and maximal submodules of $S^t/T$. From an $R$-module isomorphism $S^t/T\cong S^{T-1}$ then follows, using again Corollary \ref{maximal}, that  $1+q+\dots+q^{t-2}$ maximal submodules of $S^t$ contain fixed simple submodule $T$. Thus, $q^{t-1}$ maximal submodules do not contain $T$.   By Lemma  \ref{coatom} we see that \[
\mu_R(S^t)=-q^{t-1} \mu_R(S^{t-1}).\]
The conclusion then follows by induction. 
\end{dukaz}
\begin{Remark}
    Recall from \ref{SecDist} that a distributive semisimple module is a direct sum of pairwise non-isomorphic simple modules. In particular, for a distributive module $M$, the above calculation shows that $\mu_R(M)\in \{0,1,-1\}$. This is a special case of a general combinatorial property of distributive lattices; see [12, Example 5.1].
\end{Remark}

Using the structure of the proof of Lemma \ref{formula}, we prove that Morita equivalence preserves the Möbius function. 
\begin{Lemma}\label{Morita}
    Let $R$ and $R'$ be two Morita equivalent rings, and let $G: Mod\text-R\to Mod\text-R'$ be an equivalence of categories. Let $M\in Mod\text-R$ be a module with finitely many submodules.

    Then $G(M)$ has finitely many submodules and $\mu_{S}(G(M))=\mu_R(M)$.    
\end{Lemma}
\begin{dukaz}
    Let $S$ be a  simple $R$-module and $t$ a natural number. We first prove the statement for modules of form $S^t$. Because equivalence preserves direct limits, we get that $G(S^t)\cong G(S)^t$.  By Schur's lemma, $G(S)$ is a simple $R'$-module if and only if  $S$ is a simple $R$-module. Because equivalence is a full and faithful functor, there is a bijection $Hom_R(S, S)\leftrightarrow Hom_{R'}(G(S), G(S))$, i.e., $G$ preserves sizes of endomorphism rings. The statement then follows from Lemma \ref{formula}.
    
   Now let $M\cong S_1^{t_1}\oplus \dots S_n^{t_n}$ where  $S_i$ are simple pairwise non-isomorphic modules. Then we get
\[G(M)\cong G(S_1)^{t_1}\oplus \dots \oplus G(S_n)^{t_n},\]
    where  $G(S_i)$ are pairwise non-isomorphic simple modules. Thus, $G(M)$ is a semisimple module with finitely many submodules, and the statement then follows from Lemma \ref{productlemma} and Corollary \ref{Semisimple}.

    The fact that Morita equivalence preserves finiteness of $\mathcal{L}(M)$ then follows from Proposition \ref{Submodule-finite}.
\end{dukaz}
\subsection{Möbius inversion formula}\label{SecInv}

This brief section discusses the Möbius inversion formula. The following is a reformulation of a \textit{poset version} of the Möbius inversion formula  [12, Prop. 2] for modules.
\begin{Prop}\label{inversion}
    Let $M$ be a module with finitely many submodules. Let $f, g$ be real-valued functions on $\mathcal{L}(M)$ such that the value $g(M)$ equals the sum of values of $f$ on all submodules of $M$. Then    
\[f(M)=\sum_{\substack{N\leq M}}  g(N) \mu_R(M/N).\]
\end{Prop}

Recall that the radical of a finite-length module is zero if and only if such a module is semisimple. Thus, for a submodule $N\leq M$, factor $M/N$ is semisimple if and only if $rad~M\leq N$. Using Lemma \ref{non-semisimple}, we get the following reformulation of the Möbius inversion formula: 
\[f(M)=\sum_{\substack{rad~M\leq N\leq M}}  g(N) \mu_R(M/N).\]

\subsection{Möbius function for finite-dimensional algebras}\label{Algebras}

This subsection calculates the Möbius function for $K$-linear representations of \textit{bound quiver algebras} for acyclic quivers, in the sense defined in [1]. By Proposition \ref{thin}, if the field $K$ is infinite,  a representation has only finitely many subrepresentations if it is thin.   Note that if $K$ is algebraically closed, then any finite-dimensional $K$-algebra is Morita equivalent to some bound path algebra [1, Thm. II.3.7], and that Möbius function is preserved under Morita equivalence by Lemma \ref{Morita}.

 Let $Q$ be a finite acyclic quiver, $K$ a field,  $I$ an admissible ideal in $KQ$ and let $M\in rep_k(Q, I)$ be a nonzero representation with finitely many submodules.  

Following Lemma \ref{ASSSocle}, $\mu(M)=0$ if and only if $M$ contains a nonzero structural map. 

Now assume that $M$ is semisimple, i.e., all structural maps are zero,  and set $a_i:=dim(M_i)$. By Lemma \ref{formula}
\[\mu(M)=\prod_{\substack{1\leq i\leq n}}(-1)^{a_i} q^{\frac{a_i(a_i-1)}{2}},\]
where $q=|K|$ if $K$ is finite and $q=1$ otherwise.

\medskip 

\noindent\textbf{Aknowledgements:}
 The author wishes to express his thanks and appreciation to Jan Žemlička for his guidance and valuable insights during the research.

\medskip 

\noindent \textbf{Bibliography}

[1]~~Assem, I.   Simson, D. Skowronski, A. (2006). Elements of the Representation Theory of Associative Algebras: Techniques of Representation Theory.  \textit{LMS Student Texts}. Cambridge: Cambridge University Press. vol. 65. https://doi.org/10.1017/CBO9780511614309

\smallskip 

[2]~~Burnside, W. (2012). Theory of Groups of Finite Order. Cambridge: \textit{Cambridge University Press.}
https://doi.org/10.1017/CBO9781139237253

\smallskip 

[3]~~Camillo, V. (1975). Distributive Modules. \textit{J. Algebra}. 36(1). 16-25. https://doi.org/10.1016/0021-\\8693(75)90151-9

\smallskip 

[4]~~Enochs, E.,  Pournaki, M. R., Yassemi, S. (2024) A necessary and sufficient condition for a direct sum of modules to be distributive. \textit{Comm. Alg.} 52(2).  900-907. https://doi.org/10.1080/00927872.2023.\\2252516

\smallskip 

[5]~~Greferath, M, Schmidt, S. E. (2000). Finite-Ring Combinatorics and MacWilliams' Equivalence Theorem. \textit{J. Comb. Theory}.  Series A 92. 17-28. https://doi.org/10.1006/jcta.1999.3033

\smallskip 

[6]~~Honold, T. (2001). Characterisation of finite Frobenius rings. \textit{Arch. Math.} 76(6). 406–415. https://\\doi.org/10.1007/PL00000451

\smallskip

[7]~~Honold, T.   Nechaev, A. A. (1999), Weighted modules and representations of codes. \textit{Problems Inform. Transmission}. 35(3). 205 - 223.

\smallskip 

[8]~~Iovanov, M.  C. (2016).  Frobenius–Artin algebras and infinite linear codes. \textit{J. Pure Appl. Algebra} 220(2).  560-576. https://doi.org/10.1016/j.jpaa.2015.05.030

\smallskip 

 [9]~~Iovanov, M.  C. (2022).  On Infinite MacWilliams Rings and Minimal Injectivity Conditions. \textit{Proc. Am. Math. Soc.} 150(11).  4575-4586. https://doi.org/10.1090/proc/15929

\smallskip 

[10]~~Krasula, D. (2024). Endomorphism rings of simple modules and block decomposition. \textit{ArXiv}. https://arxiv.org/abs/2406.13582

\smallskip 

[11]~~Leinster, T.  (2012). Notions of Möbius inversion. \textit{ArXiv}.  	
https://doi.org/10.48550/arXiv.1201.\\0413

\smallskip

[12]~~Rota, G.C. (1964). On the foundations of combinatorial theory I. Theory of Möbius Functions. \textit{Probab. Theory Relat. Fields}. 340–368. https://doi.org/10.1007/BF00531932

\smallskip 

[13]~~Stephenson, W. (1974). Modules Whose Lattice of Submodules is Distributive. \textit{Proc. Lond. Math. Soc.}  (3)
28. 291–310.  https://doi.org/10.1112/plms/s3-28.2.291

\smallskip 

[14]~~Wood, J. A. (1999). Duality for Modules over Finite Rings and Applications to Coding Theory. \textit{Am. J. Math.} 121(3).  555–575. http://www.jstor.org/stable/25098937

\end{document}